\documentclass[11pt]{amsart}
\usepackage{}
\usepackage{amssymb}
\theoremstyle{plain}
\newtheorem{Thm}{Theorem}

\newtheorem{Lem}[Thm]{Lemma}

\begin{document}

%begin Topmatter
\title[Ricci expanders and type III Ricci flow]
{Ricci expanders and type III Ricci flow}

\author{Li Ma}

\address{Department of Mathematical Sciences, Tsinghua University,
 Peking 100084, P. R. China}
\email{nuslma@gmail.com}

\thanks{The research is partially supported by the National Natural Science
Foundation of China 10631020 and SRFDP 20090002110019}

\begin{abstract}
In this paper, we study how to get the Ricci expanders from
$W_+$-functional through the heat kernel estimate of the conjugate
heat equation to the type III singularity of Ricci flow. The
Gaussian upper and lower bounds are established for the related
heat kernel in accordance to the interesting work of Cao-Zhang for
the type I Rici flow.


\end{abstract}

 \maketitle

\section{Introduction}

In this paper, we study the question how to get the Ricci
expanders from $W_+$-functional through the heat kernel estimate
of the conjugate heat equation to the type III singularity of
Ricci flow (\cite{RH2}). As we shall see, the gradient estimate is
an important step in solving the Hamilton conjecture. Given a
complete non-compact Riemannian manifold $(M, g)$ of dimension n.
Recall that a family $(g(t))$ of Riemannian metrics on $M^n$ is
called a Ricci flow if $g(t)$ satisfies the following Ricci flow
equation
\begin{equation}\label{eq:ricci}
\partial_t g_{ij}(t)=-2R_{ij}(g(t)), \quad on \ M,
\end{equation}
with $g(0)=g$. We shall assume the Ricci flow $(M,g(t))$ a type
III singularity in the sense that the flow is $k$-non-collapsed on
all scales for some constant $k>0$ (see the work of G.Perelman
\cite{P02}) and is defined in $(0,\infty)$ with the curvature
bound
\begin{equation}\label{eq:bound}
|Rm(g(t))|\leq \frac{A}{A+t}
\end{equation}
for some uniform constant $A>0$. Then the conjugate heat equation
associated to the Ricci flow $(g(t))$ is
\begin{equation} \label{eq:c-heat}
 u_t=-\Delta u+Ru
\end{equation}
where $\Delta$ is the Laplacian-Beltrami operator of the evolving
metric $g(t)$. Note that the adjoint equation to the conjugate
heat equation associated to Ricci flow is
\begin{equation} \label{eq:heat}
 u_t=\Delta u.
\end{equation}
We remark that the maximum principle is true for this heat
equation associated to the Ricci flow $(M,g(t))$ (see
\cite{Shi89b}). Since we are studying the heat equation on
complete non-compact Riemannian manifold, the heat kernel is
restricted to the minimal fundamental solution $G=G(z,l; x,t)$,
$l<t$, to (\ref{eq:c-heat}) (or $u(x,t)==G(z,l;x,t)$ to
(\ref{eq:heat})). Note here that we have used an important
observation that the fundamental solution to the heat equation
(\ref{eq:heat}) can also be regards as the fundamental solution to
the conjugate heat equation (\ref{eq:c-heat}). As R.Hamilton
\cite{RH2} expected, the type III singularity of Ricci flow gives
an expanding soliton. For the $W_+$-functional
\begin{equation}\label{eq:W}
W_+(g,G,\sigma)=\int_M[\sigma(|\nabla u|^2+R)-f+n]udv_g, \; \;
\sigma=t-T>0, \; \;
\end{equation}
with
$$
G=e^{-f}/(4\pi \sigma)^{n/2},
$$
introduced by M.Freldman, T.Ilmanen, L.Ni \cite{FIN} (following
the W-functional of G.Perelman), we have
$$
\frac{d}{dt}W_+(g(t),G(t),t-T)=\int_M2(t-T)|Rc+D^2f+\frac{g}{2(t-T)}|^2dv_g\geq
0.
$$
Fix $T=0$. Hence one may get an Ricci expander by studying $W_+$
along the heat kernel and the limit of some normalization of
$g(t)$ as $t\infty$.

Our main result is below.
\begin{Thm}\label{thm:4} Let $(M,g(t))$, $t\in [0,\infty)$, be a non-flat, type
III k-non-collapsed Ricci flow with non-negative Ricci curvature
and for some $A>0$,
$$
|Rm(g(t))|\leq \frac{A}{A+t}, \ \  t>0.
$$
Then at any point $x_0\in M$, a sequence of times $\tau_k\to
\infty$, and a sequence of re-scaled metrics
$$
g_k(x,s)=\tau_k^{-1}g(x,s\tau_k)
$$
such that the pointed Ricci flow sequence $(M,x_0,g_k)$ converges
to a non-flat gradient expanding Ricci soliton in the sense of
Cheeger-Gromov sense.
\end{Thm}

This result is used in \cite{M}, where part of Hamilton's
conjecture has been proved and a full study about the
singularities of Ricci flow with positive Ricci pinching condition
is presented. Related work by us using Yamabe flow is \cite{MC}.

We now make some remarks about the proof. At the first step, we
should make sure that the $W_+$-functional is well-defined. This
will be achieved by obtaining gradient estimate of the positive
solution $u$ and the upper and lower bounds of $u$ in space-time.
According \cite{Z} (see also the work of Cao-Hamilton \cite{CH},
both following the idea of R.Hamilton), for any positive solution
$u$ to (\ref{eq:heat}) on $M\times [0,T]$, we have
\begin{equation}\label{gradient}
|\nabla\log u(x,t)|\leq
\sqrt{\frac{1}{t}}\sqrt{\log\frac{M}{u(x,t)}}
\end{equation}
for $M=sup_{M\times [0,T]} u$ and $(x,t)\in M\times [0,T]$ and
moreover, for any $\delta>0$, $t_0<t\leq T$, and $x,y\in M$:
\begin{equation}\label{eq:bound2}
u(y,t)\leq C_1u(x,t)^{1/(1+\delta)}M^{\delta/(1+\delta)}
e^{C_2d(x,y,t)^2/t}
\end{equation}
where $C_1$ and $C_2$ are positive constants depending only on
$\delta$,$d(x,y,t)$ is the distance between $x$ and $y$ in the
metric $g(t)$. We shall follow the recent interesting works of
Zhang \cite{Z}, Cao-Zhang \cite{CZ} to derive the desired heat
kernel estimates and the result is presented in section
\ref{sect2}. Our main heat kernel estimate for the Ricci flow in
Theorem \ref{thm:4} is the same result as Theorem 3.1 in
\cite{CZ}. In section \ref{sect3} we give the argument of Theorem
\ref{thm:4}.

\section{Heat kernel estimate along the Ricci flow}\label{sect2}
Let us see some analytic parts of the Type III singularity of the
Ricci flow (M,g(t)).

1.  (M,g(t)) has a space-time doubling property. Namely, the
distances of two points $x,y\in M$ at two different times $t>s>0$
are comparable in terms of $s/t$ (which is comparable in the sense
that it is bigger than a small constant). In fact, let
$\gamma(\tau;g(t))$ be the minimizing geodesic connecting $x$ and
$y$ in the metric $g(t)$. Then $d(x,y;t)=L(\gamma;g(t))$. Then by
using
$$
0\leq Rc(g(t))\leq \frac{A}{A+t},
$$
we get
$$
0\geq \frac{d}{dt} d(x,y;t)=-\int Rc(\gamma',\gamma')d\tau\geq
-\frac{A}{A+t}d(x,y;t),
$$
where the derivative is in Lipschitz sense. The latter implies
that
$$
(s/t)^{A}\leq \frac{d(x,y;s)}{d(x,y;t)}\leq 1.
$$

2. Similarly, local volume comparable property is also true. In
fact, we have
$$
0\geq \frac{d}{dt}\int_{B(x,\sqrt{t_3};t_4)}
dv_{g(t)}=-\int_{B(x,\sqrt{t_3};t_4)} Rdv_{g(t)}\geq
-\frac{A}{A+t}\int_{B(x,\sqrt{t_3};t_4)} dv_{g(t)}.
$$
Upon integration we know that the volumes of the balls
$(x,\sqrt{t_3};t_4)$ in terms of the metric $g(t_5)$ are
comparable for $t_3,t_4,t_5\in [s,t]$.

3. According to the result of E.Hebey we know that the following
Sobolev inequality holds for the k-non-collapsed $(M,g(t))$ with
$Ric(g(t))\geq 0$. Namely for all $v\in H^1(B(x,r;t))$, we have
$$
(\int |v|^{2n/(n-2)})^{(n-2)/n}\leq
\frac{c_nr^2}{|B(x,r;t)|^{2/n}}\int [|\nabla v|^2+r^{-2}
v^2]dv_{g(t)}.
$$

We shall choose $r=c\sqrt{t}$ for $c\in (0,1)$. Then by assumption
$|Ric|\leq \frac{A}{A+t} $ and the k-non-collapsing property, we
know that
$$
|B(x,r;t)|\geq kA^{-n}t^{n/2}.
$$
Hence, we obtain the following uniform Sobolev inequality along
the Ricci flow. For all $v\in H^1(B(x,\sqrt{t};t))$, we have
$$
(\int |v|^{2n/(n-2)})^{(n-2)/n}\leq \frac{c_nA^2}{k^{2/n}}\int
[|\nabla v|^2+t^{-1} v^2]dv_{g(t)}.
$$

In the following we try to get the Gaussian upper and lower bounds
for the heat kernel to the heat equation. The general program for
this work is in three steps below. 1). We derive a weaker
on-diagonal upper bound
$$
u(x,t)\leq \frac{const.}{f(t)}
$$
for some increasing function in t. Here $u(x,t)=G(x_0,0;x,t)$. The
method is the Moser iteration and in this step, the uniform
Sobolev inequality plays a important role. 2). We derive the
Gaussian upper bound by the exponential weight method due to
E.B.Davies. 3). We derive a on-diagonal lower bound at some point
and then we obtain the full Gaussian lower  bound by using the
gradient estimate obtained by Zhang \cite{Z} and Cao-Hamilton
\cite{CH}, both papers follow the idea of the work of R.Hamilton
for standard heat equation.

However, such a machinery can be go through in our setting. So we
shall prove the Gaussian upper and lower bound of the heat kernel
in a new way (which is also observed by X.Cao and Zhang
\cite{CZ}).

Once we have the Gaussian upper and lower bound of the heat
kernel, we immediately see the $W_+$-functional is well-defined
along the Ricci flow.

Let us now go to the detail. Note that if $G=G(z,l;x,t)$, $l<t$,
is the fundamental solution of
$$
u_l=-\Delta u+Ru
$$
along the Ricci flow, then as a function of $(x,t)$,
$p(x,t;x_0,l):=G(x_0,l;x,t)$, $l<t$, is the fundamental solution
to the heat equation
$$
 u_t=\Delta u.
$$

Let $u=u(x,t)$ be a positive solution to the heat equation
(\ref{eq:heat}) in the region
$$
Q_{\sigma \tau}:=\{(y,s)\in M\times [\tau-(\sigma r)^2,\tau,
d(y,x;s)\leq \sigma r]\}.
$$
Here $r=\sqrt{t}/8>0$, $\sigma\in [1,2]$. Then for any $p\geq 1$
we have
$$
\partial_t u^p\leq \Delta u^p.
$$

Choose a non-negative smooth function $\phi:[0,\infty)\to [0,1]$
such that
$$
|\phi'|\leq \frac{2}{(\sigma-1)r}, \ \ \phi'\leq 0,
$$
and $\phi(\rho)=1$ for $0\leq \rho\leq r$ and $\phi(\rho)=0$ for
$\rho\geq \sigma r$. Then we choose a smooth non-negative function
$\eta$ such that
$$
|\eta'|\leq \frac{2}{(\sigma-1)^2r^2}, \ \ \eta'\geq 0,
$$
and $\eta(\rho)=1$ for $\tau-r^2\leq s\leq \tau$ and $\phi(s)=0$
for $s\geq \tau-(\sigma r)^2$.

  Define $\xi=\phi(d(x,y;s))\eta(s)$.

Set $w=u^p$ and using $w\xi^2$ as a testing function to the above
differential inequality we deduce that
$$
\int \nabla w\cdot \nabla (w\xi^2)dv_{g(s)} ds\leq -\int
\partial_s w w\xi^2 dv_{g(s)} ds.
$$

Note that the left hand side is
$$
\int |\nabla (w\xi)|^2dv_{g(s)}ds-\int |\xi|^2w^2dv_{g(s)}ds
$$
and the right hand side is
$$
\int w^2\xi^2\partial_sdv_{g(y,s)}ds-\frac{1}{2}\int (w\xi)^2
Rdv_{g(s)}ds -\frac{1}{2}\int (w\xi)^2dv_{g(y,\tau)},
$$
Using $R\geq 0$ and our choice of $\phi$ and $\eta$ we know that
the latter is bounded by
$$
\frac{c}{(\sigma-1)^2 r^2}\int w^2dv_{g(s)}ds-\frac{1}{2}\int
(w\xi)^2dv_{g(y,\tau)}.
$$

Re-arranging the above relations we get
$$
\int |\nabla (w\xi)|^2dv_{g(s)}ds+\frac{1}{2}\int
(w\xi)^2dv_{g(y,\tau)}\leq \frac{c}{(\sigma-1)^2 r^2}\int
w^2dv_{g(s)}ds.
$$
Note that the Sobolev inequality gives us
$$
(\int |w\xi|^{2n/(n-2)})^{(n-2)/n}\leq \frac{c_nA^2}{k^{2/n}}\int
[|\nabla (w\xi)|^2+r^{-2} (w\xi)^2]dv_{g(t)}.
$$
by Holder inequality
$$
\int (w\xi)^{2(1+2/n)}dv_{g(s)}\leq \int
(w\xi)^{2n/(n-2))}dv_{g(s)})^{(n-2)/n}(\int
(w\xi)^{2}dv_{g(s)})^{2/n}.
$$
Then using the trick as in \cite{DWY} we get that
$$
\int_{Q_r(x,\tau)} w^{2\theta}\leq
c(k,A)(\frac{1}{(\sigma-1)^2r^2}\int_{Q_{\sigma
r}(x,\tau)}w^2)^{\theta}
$$
with $\theta=1+2/n$. Choose
$$
\sigma_0=2, \sigma_i=2-\sum_{1}^i 2-{j}, \ \ p=\theta^i,
$$
and we find that
$$
\sup_{Q_{r/2}(x,\tau)}u^2\leq
\frac{c(k,A)}{r^{n+2}}\int_{Q_r(x,\tau)} u^2dv_{g(s)} ds.
$$
Using the general trick of Li-Schoen \cite{LS}, we then find the
$L^1$ mean value inequality in the form below:
$$
\sup_{Q_{r/2}(x,\tau)}u\leq
\frac{c(k,A)}{r^{n+2}}\int_{Q_r(x,\tau)} udv_{g(s)} ds.
$$
Using $u(x,t)=G(x_0,0;x,t)$, $r=\sqrt{t}$ and the fact $\int_M
udv_{g(s)}=1$, we obtain
$$
G(x_0,0;x,t)\leq  \frac{c(k,A)}{t^{n/2}}.
$$

Next we prove the lower bound by the trick of Perelman \cite{P02}.
This argument is borrowed from \cite{CZ}.  Let
$u=u(x,t)=G(x,t;x_0,t_0)$ for $t<t_0$. Claim that for some uniform
constant $C>0$, we have, for $t<t_0$,
$$
G(x_0,t;x_0,t_0)\geq
\frac{C}{\tau^{n/2}}e^{-\frac{1}{2\sqrt{\tau}}\int_t^{t_0}
\sqrt{t_0-t}R(x_0,s)ds}
$$
where $\tau:=t_0-t$.

Following G.Perelman, we set
$$
u=(4\pi \tau)^{-n/2}e^{-f}.
$$
 Using Perelman's differential Harnack inequality for the fundamental
 solution we have that for $\gamma(t)=x_0$, we have
$$
-\partial_tf(x_0,t)\leq
\frac{1}{2}R(x_0,t)-\frac{1}{2\tau}f(x_0,t).
$$
Then for any $t_2<t<t_0$ we can integrate the above inequality to
obtain
$$
f(x_0,t_0)\sqrt{t_0-t_2}\leq
f(x_0,t_1)\sqrt{t_0-t_1}+\frac{1}{2}\int_{t_2}^{t_1}\sqrt{t_0-s}R(x_0,s)ds.
$$
We remark that by the asymptotic formula for $G$ we know that for
$t_1$ approaches to $t_0$, $f(x_0,t_1)$ stays bounded since
$G(x_0,t_1;x_0,t_0)(t_0-t_1)^n/2)$ is bounded between two positive
constants. Then for any $t\leq t_0$,
$$
f(x_0,t)\leq \frac{1}{2\sqrt{t_0-t}}\int \sqrt{t_0-s}R(x_0,s)ds.
$$

Hence,
$$
G(x_0,t;x_0,t_0)\geq c (4\pi
\tau)^{-n/2}e^{-\frac{1}{2\sqrt{t_0-t}}\int
\sqrt{t_0-s}R(x_0,s)ds}.
$$
Using the assumption $|R(x,s)|\leq A/(t_0-s)$, we then know that
$$
G(x_0,t;x_0,t_0)\geq c (4\pi \tau)^{-n/2}.
$$

In summarize, we have obtained the below.

\begin{Lem}\label{keylemma} Let $(M,g(t))$,$t\in [0,\infty)$,
be a k-non-collapsed Ricci flow with bounded curvature $|Rm|\leq
\frac{A}{A+t}$, and non-negative Ricci curvature. Then there exist
positive constant $C_1$ and $c_2$ which depends only on $k, n, A$,
such that for all $x,x_0\in M$, $t>0$, we have
$$
G(x_0,0;x,t)\leq \frac{C_1}{t^{n/2}},
$$
and
$$
G(x_0,0;x_0,t)\geq \frac{C_2}{t^{n/2}},
$$
\end{Lem}

With this Lemma to replace Theorem 2.1 in \cite{CZ}, we can get
the same result as Theorem 3.1 in \cite{CZ}, which is enough for
our use of $W_+$-functional.

\section{Blow-down for the Type III Ricci flow}\label{sect3}

Choose suitable time sequence $\tau_k\to\infty$ and point sequence
$x_k\in M$ for the blowing up metrics as in Hamilton \cite{RH2}.
Consider the pointed Ricci flow $(M,g_k,x_k)$ with
$$
g_k(s):=\tau_k^{-1}g(\cdot,t_k+s\tau_k).
$$
Let, for $s\in [1,4]$ and for $x_0=x_k$,
$$
u_k=u_k(x,s):=\tau_k^{n/2}G(x,s\tau_k;x_0,6\tau_k).
$$
Then $u_k$ satisfies
$$
\partial_su_k=-\Delta_{g_k}u_k+R(g_k)u_k.
$$
Recall that $f_k$ is defined by the relation $$ (4\pi
s)^{-n/2}e^{-f_k}=u_k.
$$ Using the upper bound for $u_k$, we know that
$$
-f_k=\log u_k+\frac{n}{2}\log(4\pi s)\leq C_0
$$
for all $k=1,2,...$ and $s\in [1,3]$. Since $\int_Mu_kdg_k=1$, by
our uniform bounds for $u$, we know that there is a limit
$u_{\infty}$ of $u_k$ as $k\to\infty$. Note also that
$$
W_{+k}(s)\leq C-n
$$
for all $k=1,2,...$ and $s\in [1,3]$. By scaling we know that
$$
W_{+k}(s)=W_+(g,u,s\tau_k)\leq C-n,
$$
where $u=u(x,t)=G(x,t;x_0,6\tau_k)$, $0\leq t\leq 4\tau_k$. Note
that using the asymptotic behavior of $u$ (see \cite{Chow}) we
have
$$
W_+(g,G(x,t; x_0,6\tau_k),s\tau_k)\leq
W_+(g,G(x,t;x_0,6\tau_{k+1}),s\tau_k)+\circ_k(1),
$$
and using the increasing property we have
$$
W_+(g,G(x,t;x_0,6\tau_{k+1}),s\tau_k)\leq
W_+(g,G(x,t;x_0,6\tau_{k+1}),s\tau_{k+1}).
$$
 Hence there exists the limit $W_{\infty}(s)$ for
the sequence ${W_{+k}}(s)$ as $k\to\infty$.

Then we can get an expanding Ricci soliton similar to Zhang did in
the case for shrinking soliton \cite{Z}. Since the argument is
almost the same, we omit the detailed proof. Thus we have
completed the proof of Theorem \ref{thm:4}.

\emph{Acknowledgement}: The author would also like to thank IHES,
France for host and the K.C.Wong foundation for support in 2010.

\end{document}